\documentclass[oumk,xs]{xxllp}
\usepackage{amsmath,amssymb, amsfonts} 
\usepackage{url}

\usepackage{xxllp}
\usepackage{LLPthm}
\usepackage{xcolor}

\newtheorem{theorem}{Theorem}[section]

\theoremstyle{remark}

\theoremstyle{definition}
\newtheorem{defn}{Definition}[section]

\newcommand\JLLuk{\mathsf{J}_{\Lslash}}
\newcommand{\Lslash}{\mathchoice
{\mbox{\sf\L}}
{\mbox{\sf\L}}
{\mbox{\sf\scriptsize\L}}
{\mbox{\sf\tiny\L}}}

\usepackage{apalike}

\begin{document}

\twoAuthorsTitle{Nancy Abigail Nuñez Hernández}{Francisco Hernández Quiroz}{Justification logic and the epistemic contribution of deduction}

\allowdisplaybreaks

\noindent \textbf{Abstract:} Accounting for the epistemic contribution of deduction has been a pervasive problem for logicians interested in deduction, such as, among others, Jakko Hintikka. The problem arises because the conclusion validly deduced from a set of premises is said to be ``contained'' in that set; because of this containment relation, the conclusion would be known from the moment the premises are known. Assuming this, it is problematic to explain how we can gain knowledge by deducing a logical consequence implied in a set of known premises. To address this problem, we offer an alternative account of the epistemic contribution of deduction as the process required to deduce a conclusion or a theorem, understanding such process not only in terms of the number of steps in the derivation but, more importantly, the reason for or justification for every step. That is, we do not know a proposition unless we have a justification or proof to hold that proposition. With this goal in mind, we develop a justification logic system which exhibits the epistemic contribution of a deductive derivation as the resulting justified formula.




\section{Introduction}

In 1970, Jaakko Hintikka referred to the problem of accounting for the epistemic contribution of deduction as the ``scandal of deduction.''\footnote{The problem had been acknowledged by others before Hintikka: \cite[Book II, chap. iii, \S 1]{Mill1843}, \cite{Cohen&Nagel} who referred to it as ``the paradox of inference,'' or \cite{Dummett1978} who labeled it as ``the justification of deduction.'' According to \cite{parikh2008}, Plato had already alluded to this problem in \textit{Meno}. Here we take Hintikka's exposition of the problem as our starting point because it is the one that has received the most attention.} According to Hintikka, the problem seems to arise from the ``tautological'' or ``analytical''\footnote{Because of the deduction theorem, whenever a set of premises $\Gamma$ entail a conclusion $\alpha$, the conditional $\Gamma \to \alpha$ is a tautology; thereby, we have the tautological nature of deduction. To understand how such a ``tautological'' nature of deduction has been equated to its presumed ``analytical'' nature, see \cite[Chapter~2]{primiero2007information}.} nature of deduction: the conclusion validly deduced from a set of premises is said to be ``contained'' in that set. Because of this relation of containment, it has been problematic to explain how deducing the conclusion implied by a set of known premises yields an epistemic contribution. 

However, even if it is problematic to explain the epistemic contribution, it is actually the case that deduction yields an epistemic contribution for humans and other non-idealized epistemic agents bounded by the limits of computational complexity. Given these limits, humans are not logically omniscient, that is, we do not know all the logical consequences implied by our current knowledge, but we can actually get to know some them if we deduce them. Deducing such logical consequences leads to know them not only because it allows us to be aware of them, but also because it provides us with the reasons or justification to accept them. In view of these facts, the aim of this paper is to offer an alternative account of the epistemic contribution of deduction in terms of the process required to deduce a specific conclusion or a theorem, understood not only in terms of the number of steps in the derivation, but also as the reasoning to justify every step. To explicitly exhibit such process, in the third section of this paper, we present a justification logic system which will allow us to indicate the epistemic contribution of a deductive derivation as the resulting justified formula.

The framework of justification logic originated in the logic of proofs $\mathsf{LP}$ developed by Sergei Artemov \cite{artemov1994,artemov1995}; the rough idea behind $\mathsf{LP}$ is to ``extend classical propositional logic by adding symbolically represented proofs into the language of the system'' \cite[1061]{artemov2005introducing}. $\mathsf{LP}$ is inspired by constructive (intuitionistic) mathematics, where truth is understood in terms of the existence of a proof, that is, a mathematical statement is accepted as true if it has a proof. In epistemic terms, this means that it is necessary --though not sufficient-- to have reasons, a proof or a justification to know something. $\mathsf{LP}$ intends to capture this idea and express it in the language: the existence of a proof for a statement $S$ that is true would be formalized as $t : S$, thereby expressing within the language that $t$ is a proof for $S$. 

Thus, $\mathsf{LP}$ ``may be regarded as a device that makes reasoning about knowledge explicit and keeps track of the justification'' \cite[1061]{artemov2005introducing}. Therefore, justification logic can show that the epistemic contribution of finding the justified formula $t : \varphi$ is greater than just having the formula $\varphi$ without a justification, which can be acknowledged upon noting that justified formulas of the form $t : \varphi$ are meant to convey that $t$ is a reason ---or justification--- for $\varphi$. Taking advantage of this feature in the language of justification logic, we shall present a justification logic system with enough power to deduce or derive tautologies or theorems of propositional logic, explicitly exhibiting the epistemic contribution of the derivations. Although our system has enough power to deduce all tautologies of propositional logic, it doesn't include as axioms of the system, which prevents logical omniscience. By preventing logical omniscience, we give space for the process that generates knowledge by deduction. 



The remaining parts of this paper are organized in the following way. The second section motivates our proposal by discussing some literature on the problem of accounting for the epistemic contribution of deduction. Since it is not the aim of this paper to provide a comprehensive survey of the literature on this topic, we focus on some relevant and representative reactions to the problem of accounting for the epistemic contribution of deduction; these reactions range from the philosophical (represented by Cohen and Nagel) to the technical, which include the tradition of logical responses (initiated by Hintikka), but have evolved to incorporate theoretical tools and insights from computer science and mathematics (Du{\v{z}}{\'\i}, Abramsky, D'Agostino and Floridi). Given the lack of consensus on how to account for the epistemic contribution of deduction, in the third section we provide an alternative proposal in the form of a justification logic system which allows us to deduce or derive the tautologies of propositional calculus accompanied by their respective justification, thereby rendering explicit their epistemic contribution within the language of the system.


\section{A selection of the literature}
The precise epistemic contribution of deduction has been controversial and puzzling within philosophy and scientific disciplines interested in deductive inferences and knowledge, such as cognitive psychology and other cognitive sciences. According to the standard understanding of deduction, a conclusion that is logically implied by a set of premises is ``contained'' in them. Because of this relation of containment, it has been problematic to explain how knowledge can be gained by deduction. Ian Rumfitt provides a concise description of the problem: ``some philosophers have been puzzled by how a thinker can gain knowledge by deduction. There are two problems here: first, how exercising a deductive capacity can yield knowledge; second, how it can yield knowledge that the thinker does not posses already?'' \cite[62]{Rumfitt2008}. The aim of this paper is to address the latter question ---that which seeks for an inquiry into the epistemic contribution of deduction.

Even though it seems to be widely accepted that deduction yields an epistemic contribution by allowing us to gain knowledge, it is still hard to understand exactly how this can happen if the conclusion of a valid deduction is contained in the previously known premises. As Rohit Parikh points out, the question can be traced as far back as Plato: ``In deductive reasoning, if $\phi$ is deduced from some set $\Gamma$, then $\phi$ is already implicit in $\Gamma$. But then how do we learn anything from deduction? That we do not learn anything is the (unsatisfying) answer suggested by Socrates in Plato's \textit{Meno}.'' \cite[459]{parikh2008}. However, as Parikh notes, such an answer is not only unsatisfying, but also untenable given its commitment to logical omniscience.\footnote{It is not among the goals of this paper to discuss the problem of logical omniscience. Suffice it to say that philosophers and logicians working on epistemic logic understand logical omniscience as an agent’s ability to know all the logical consequences that follow from her knowledge (and to know all the tautologies as well); since no real agent possesses such ability, logical omniscience is considered a problem.}

The epistemic contribution of deduction has been considered so problematic that Morris Raphael Cohen and Ernest Nagel considered that it led to a paradox: 
\begin{quote}
    If in an inference the conclusion is not contained in the premises, it cannot be valid; and if the conclusion is not different from the premises, it is useless; but the conclusion cannot be contained in the premises and also possess novelty; hence inferences cannot be both valid and useful \cite[p.~173]{Cohen&Nagel}.
\end{quote}
Cohen and Nagel's way of dealing with the paradox is to point out that the sense of novelty we experience in deducing something is merely psychological \cite[p.~174]{Cohen&Nagel}; the agents may be surprised by the conclusion, but this does not mean that its deduction yields any epistemic contribution.

Jaakko Hintikka criticizes the kind of response provided by Cohen and Nagel to the paradox, writing that
\begin{quote}
    [i]f no objective, non-psychological increase of information takes place in deduction, all that is involved is merely psychological conditioning, some sort of intellectual psycho-analysis, calculated to bring us to see better and without inhibitions what objectively speaking is already before our eyes \dots\ [and] all the multifarious activities of a contemporary logician or mathematician that hinge on deductive inference are as many therapeutic exercises calculated to ease the psychological blocks and mental cramps that initially prevented us from being, in the words of one of these candid positivists, ``aware of all that we implicitly asserted'' already in the premises of the deductive inference in question \cite[p.~135]{hintikka1970}.
\end{quote}
Hintikka acknowledges that deduction yields an epistemic contribution, and he  refers to the long-standing failure to account for it as the scandal of deduction, pointing out that we cannot even properly answer first year students who ask how we could learn from deduction, given that it is ``{}`tautological' or `analytical' and that logical truths have no `empirical content'.'' \cite[p.~135]{hintikka1970} Hintikka realizes that if deduction only ``brings us to see what is already before our eyes'' then it is hard to make sense of why it often seems that we actually do not know these things. 

To address the problem, Hintikka develops a proposal in which he distinguishes between depth information and surface information; deduction can only increase the latter in a set of deductions within polyadic predicate calculus, but neither within monadic predicate calculus (a fragment of first-order logic in which all relation symbols take only one argument and there are no function symbols) nor propositional logic.\footnote{It is worth mentioning that Hintikka's proposal is framed under the classical semantic information theory developed by Yehoshua Bar-Hillel and Rudolf Carnap \cite{Carnap1952CSI}. According to this theory, logical truths are uninformative while contradictions are maximally informative.

The reason underlying this apparently counterintuitive result goes back to Carnap's notion of intension. The intension of a declarative sentence is identified as the set of possible worlds in which the sentence is true; tautologies are true in every possible world and contradictions are false in every possible world. The semantic information contained in a sentences is identified as the set of possible worlds in which the sentence is false, which are the worlds excluded by the truth of the sentence. Accordingly, the semantic information contained in a sentence $s$ is expressed as its content, which is denoted as `cont' and is defined as the measure $m$ (which assigns the same probability to every world) of the complement of the sentence $s$:
\begin{center} 
{cont($s$)} = $1-m(s)$
\end{center}
A tautology $t$ is a sentence that is logically true in every possible world: 
\begin{center}
$m(t)=1$
\end{center}
So, according to the Bar-Hiller-Carnap theory of semantic information, the content of a tautology $t$ (the number of possible worlds it excludes) would be:
\begin{center} 
cont$(t)=1-1=0$
\end{center} 
Hintikka realized that if tautologies are uninformative, deduction is also uninformative because in classical logic a sentence $\psi$ can be deduced from a finite set of premises $\phi_1, \dots, \phi_n$ if and only if the conditional $\phi_1, \dots, \phi_n \to \psi$ is a tautology. Under the Bar-Hiller-Carnap theory of semantic information, deduction can hardly provide anything new in terms of knowledge.} This proposal has not been further developed because it fails to apply to deductions in propositional calculus.\footnote{For an extensive and detailed explanation of Hintikka's proposal, we recommend \cite{sequoiah2008}.} However, since many deductions that take place in propositional calculus yield an epistemic contribution, in the following section of this paper we shall present a justification logic system capable of recovering all tautologies from propositional calculus exhibiting the work involved in each derivation, which in turn exhibits its epistemic contribution.   

In a similar vein to Hintikka's reaction, Marie Du{\v{z}}{\'\i} reacts to the paradox of inference and the scandal of deduction asserting that ``as many fellow logicians and mathematicians will no doubt agree, the conclusion of a valid argument is often very useful\dots\ It seems evident that there is \textit{something} that we \textit{learn} when deducing the conclusion of a sound argument'' \cite[474]{Duvzi2010paradox}. To address the problem posed by the paradox of inference, Du{\v{z}}{\'\i} develops an account of information to show how deduction can provide novel and useful information. For this purpose, she appeals to a framework of procedural semantics, in which the meaning of a sentence is modeled in terms of \textit{constructions} (which are abstract, algorithmically structured procedures), which in turn are the vehicles of information. In this way, ``although analytically true sentences provide no \textit{empirical information} about the state of the world, they convey \textit{analytic information}, in the shape of constructions prescribing how to arrive at the truths in question'' \cite[473]{Duvzi2010paradox}. In this way, the novel and useful contribution of valid deduction can be acknowledged in the metalanguage: ``A \textit{construction} of the conclusion may not occur in the premises; if it does not we have to discover it, and the construction we discover is new to us, hence epistemically useful and non-trivial.'' Since such constructions are not syntactic objects, their discovery does not take place within the language in which the deduction takes place. In contrast, the proposal we present in the next section explicitly exhibits the epistemic contribution of deduction within the language.             

Along with Du{\v{z}}{\'\i}, another computer scientist puzzled by this problem is Samson Abramsky, who considers that if deduction does not yield any contribution, then it is not clear why we have to compute the logical consequences implied by a set of premises in order to know such consequences. According to him, the natural answer to the question `Why do we compute?' is ``to gain information (which we did not previously have!) But how is this possible?\dots\ Isn't the output \textit{implied} by the input?\dots\ If we extract logical consequences of axioms, then surely the answer was already there implicitly in the axioms; what has been added by the derivation?'' \cite[p.~483]{ABRAMSKY2008}. Abramsky's proposal to address this problem is framed within the field of information dynamics in computer science, and it is focused on information flow in computation, drawing upon game-based models of interactive processes.     

Although Abramsky talks about information rather than knowledge, it is easy to see the epistemological dimension in his question: following Du{\v{z}}{\'\i}, we just need to realize that if a valid argument is considered uninformative because its conclusion is contained in the premises, and uninformative arguments are considered epistemically useless, then valid arguments are epistemically useless \cite[482]{Duvzi2010paradox}. However, according to Johan van Benthem and Maricarmen Martinez, ``deduction is useful for the purpose of `extracting information' from the data at our disposal'' \cite[227]{BenthenMartinez}. Furthermore it is safe to assume that the information extracted in this way by an agent will become part of her body of knowledge, which shows that deduction can be informative and epistemically useful. Even though at the end they suggest that some work on Kolmogorov complexity or the logic of proofs includes elements to address the problem of accounting for the epistemic contribution of deduction (although they do note that the increase in information from the point of view of Kolmogorov complexity would be too short), van Benthem and Martinez's main goal is rather to discuss the role of information in logic throughout various logical developments.

Marcello D'Agostino and Luciano Floridi's \cite{floridi:dagostino:scandal} approach to the scandal of deduction deserves special attention: they propose a solution based on complexity measures, focusing specifically on the complexity of the Boolean satisfiability problem (SAT). D'Agostino and Floridi develop a hierarchy of propositional logics that are tractable; they commence by distinguishing uninformative inferences from  informative inferences. While in the former the information carried by the conclusion is obviously contained in the premises, in the latter the information carried by the conclusion is not obviously contained in the premises, and virtual information (the introduction of hypothesis or assumptions) is necessary in order to reach the conclusion. Such virtual information is not even implicitly contained in the premises of informative inferences, so every time a new assumption is introduced, information is increased. According to D'Agostino and Floridi, this would be the source of the intractability of full propositional logic.

To tackle the problem, D'Agostino and Floridi propose a family of restricted logical deduction systems. The first step towards achieving this family of deduction systems is to eliminate virtual information by eliminating discharge rules, thereby leaving out informative inferences. Although informative inferences are those which increase information, they are left out because this sense of information is ``overloaded'' and leads to logical omniscience. By eliminating discharge rules, D'Agostino and Floridi achieve a tractable system of introduction and elimination rules for propositional connectives (``intelim'' for short) that are ``analytical'' (i.e. only lead to uninformative inferences).      

Once they have achieved such a tractable system, D'Agostino and Floridi redefine the meaning of logical operators in informational terms. This new meaning given to the logical operator does not require any use of virtual information. However, by redefining the meaning of logical operators and proposing a tractable system in which problems such as satisfiability (SAT) would be tractable, D'Agostino and Floridi no longer rely on the intractability of SAT to show that deduction is informative. Moreover, their tractable system does not have the expressive power of classical propositional logic. Therefore, in the final sections of their paper, they address 
\begin{quote}
    the problem of gradually retrieving the full deductive power of classical propositional logic by means of a bounded recursive use of virtual information\dots\ its iterated incremental use leads to more and more powerful deductive systems\dots\ Although these systems are all tractable, their growing computational complexity approaches intractability as their deductive power approaches that of classical propositional logic \cite[p.~276]{floridi:dagostino:scandal}.
\end{quote}
Thus, it seems that the D'Agostino and Floridi proposal backfires on their original aim: to tackle the Scandal of Deduction accounting for the informativeness of classical propositional logic.\footnote{Moreover, in D'Agostino and Floridi's tractable system, the problem known as Logical Entailment (LE) is tractable: ``Intelim deducibility is a tractable problem. Whether a formula P is intelim deducible from a finite set $\Gamma$ of formulas can be decided in polynomial (quadratic) time'' \cite[p.~296]{floridi:dagostino:scandal}. However, in classical propositional logic, LE or deciding whether some formula is the logical consequence of a set of formulas is not a tractable problem; LE is equivalent to TAUT, which is coNP-complete \cite{arora_barak_2009}.} The proposal we develop in the next section does not face this sort of inconvenience because it is developed within the language of classical propositional logic from the beginning.  

As we can see, the problem posed by the so-called scandal of deduction and the paradox of inference has been understood in different ways and thereby has motivated the most diverse reactions. Philosophers like Cohen and Nagel understood the problem in terms of the supposed impossibility of accounting for the usefulness of valid deductions, which led them to conclude that any novelty in deduction is due to a mere psychological feeling of surprise. Unlike Cohen and Nagel, we claim that deduction can increase knowledge, and in this way, yield a significant epistemic contribution; to underpin this claim, we present a justification logic system that exhibits the epistemic contribution of derivations in the following section. 

Hintikka understood the problem in terms of accounting for the information gains through deduction, so his reaction differed from that of Cohen and Nagel. Hintikka acknowledged that deduction can give us new information, and he tried to explain this by the aforementioned distinction between surface information and depth information. This distinction has inspired many developments in logic that in turn distinguish between implicit and explicit attitudes within a logical system, for instance, between implicit and explicit knowledge. The main difference between them is that the former is closed under logical consequence, while the latter is not. Under this kind of logical framework, agents implicitly know all the logical consequences of any known proposition, so even if they derive new knowledge by deduction at the explicit level, they remain logically omniscient at the implicit level. Considering that logical omniscience is a problem even if it remains in the domain of implicit knowledge, we shall appeal to a logical framework that avoids any commitment to logical omniscience, i. e., justification logic \cite{ArtemovKuznets2009}.      

Others, such as Du{\v{z}}{\'\i}, Abramsky, or van Benthem and Martinez, have understood the problem in terms of information, so their reactions have been proposals also developed in terms of information. In contrast, the aim of this paper is to offer an alternative account in terms of knowledge and the epistemic contribution of deduction. To this end, we shall present a logical system which shows how the process of deducing logical consequences yields an epistemic contribution.

\section{The epistemic contribution of deduction in $\JLLuk$}

There are infinitely many logical consequences implied by the current knowledge of a non-ideal epistemic agent with bounded cognitive resources --no matter whether the agent is human or artificial. Although this kind of agent cannot know each and every single one of these consequences, she can known some of them if she performs the relevant deduction which, besides delivering an unknown proposition, provides reasons that justify the new piece of knowledge. In this way, deduction yields an epistemic contribution because it presents the agent an unknown proposition as well as the reasons that justify it. However, this may not be fully clear when it is considered in a system like classical propositional calculus. For example, the epistemic contribution of deriving $r$ from premises $p$, $p \to q$ and $(p \to q) \to (q \to r)$ would be contested because $r$ is ``contained'' in the premises as the consequent of the consequent of the third premise: 
\begin{align*}
  & 1.\ p   & & \text{premise}\\
  & 2.\ p \to q & & \text{premise}\\
  & 3.\ (p \to q) \to (q \to r) & & \text{premise}\\
  & 4.\ q & & \text{MP 2, 1}\\
  & 5.\ \underline{q \to r} & & \text{MP 3, 2}\\
  & 6.\ r & & \text{MP 5, 4}
\end{align*}

To address this potential objection to the epistemic contribution of deduction we developed the system $\JLLuk$ to show the epistemic contribution of deduction as the process required to deduce a conclusion or a theorem, understanding such process not only in terms of the number of steps in the derivation, but more importantly, the reason or justification for every step. The system $\JLLuk$ is based on justification logic and it has the advantage of exhibiting the epistemic contribution of a deductive derivation as the resulting justified formula.

\section{Justification Logic}
\subsection{A system to express the epistemic contribution of deduction}

The main feature of justification logics is the introduction of justification terms within the language.  Justification terms intend to \textit{explicitly} capture or express the existence of justifications or reasons to hold propositions. Thus, instead of simply asserting a proposition that is known to be true, the justification or reason for the truth of the proposition should also be asserted. For this purpose, an epistemic operator is introduced in the from of \textit{justification terms} for the propositions that are actually justified: $t: P$, which can be read as ``\textit{P is justified by reason t}'' or ``\textit{t justifies P}'', where $t$ is the justification term. Hence, working with a justification logic system will allow us to explicitly express the epistemic contribution of deduction as the resulting justified formulas from a derivation or deduction. As an illustration, let us contrast an example from classical propositional logic with one from justification logic. While in classical propositional logic it might not be clear the epistemic contribution of deducing $q$ from $p$ and $p \to q$, justification logic can show that reasons to hold the derived conclusion are not the same as the reasons to hold the premises. Let's say that $p$ and $p \to q$ are justified by $t_1$ and $t_2$ respectively, then $q$ would be justified by $t_1 \cdot t_2$.    

\begin{align*}
    & p  && t_1:\ p\\
    & \underline{p \to q}  && \underline{t_2:\ p \to q}\\
    & q  &&  t_1 \cdot t_2:\ q\\
\end{align*}

While more conventional versions of epistemic logic systems provide the formal tools to model \textit{what} the agent knows when she derives a conclusion, justification logic provides the formal tools to model \textit{why} the agent knows such conclusion, which reflects the process required for knowledge acquisition through deduction. Let us contrast justification logic with a standard system of epistemic logic $\mathsf{S4}$, in which the modal operator $K$ is used to express that an agent knows a formula $\phi$ as $K\phi$. In a justification logic system $\mathsf{J4}$ if $\phi$ is known, then it is necessary to express the reason that justifies the knowledge of $\phi$ in terms of $t: \phi$.\footnote{It is a common practice to name modal logic systems by stringing axiom names after $\mathsf{K}$; e. g. $\mathsf{KT}$, $\mathsf{K4}$, and so on, with $\mathsf{K}$ itself as the simplest case. In the names of the justification logic counterparts for such modal logic there is a substitution of $\mathsf{J}$ for $\mathsf{K}$. $\mathsf{S4}$ sometimes is called $\mathsf{KT4}$.} In this way the use of justification terms allows us to express the process or the reasons that leads to knowledge of the formula. 

However, $t: \phi$ is still pretty similar to $K\phi$. As a matter of fact, the system $\mathsf{J4}$ is analogous to $\mathsf{S4}$. The logic of proofs $\mathsf{LP}$, from which $\mathsf{J4}$ evolved, incorporates modal logic $\mathsf{S4}$ and intutionistic logic.  The modal logic axiom schemes Factivity and Positive Introspection present in $\mathsf{S4}$ are also part of $\mathsf{LP}$.\footnote{In modal logic, Factivity is represented by the axiom scheme $\square \varphi \to \varphi$, which in justification logic is $t:\varphi \to \varphi$. In modal logic, Positive Introspection is represented by the axiom scheme $\square \varphi \to \square \square \varphi$; in $\mathsf{LP}$ Positive Introspection is $t:\varphi \to !t:t:\varphi$. ! is a one-place function symbol on justification terms known as \textit{Fact Checker}; if $t$ is a justification of something, $!t$ is a justification that $t$ is such justification.} Hence, the justification logic developed by Artemov is a counterpart of the modal logic $\mathsf{S4}$. Like other modal logic systems, $\mathsf{S4}$ has the inference rule known as Rule of Necessitation: 
\begin{center}
    \underline{$\vdash \varphi$}\\
    $\vdash \square \varphi$
    
\end{center}
In view of this rule, if a formula is a tautology in the system, it is known by an epistemic agent. 
So according to the Rule of Necessitation, agents are logically omniscient because they know all the tautologies in $\mathsf{S4}$. 
Given that $\mathsf{LP}$ is equivalent to $\mathsf{S4}$ and includes all the tautologies of propositional logic, it also lead to epistemically omniscient agents who know all the tautologies of propositional logic.

In contrast, the system we present below as $\JLLuk$ is not so strong. $\JLLuk$ does not include all tautologies of propositional logic nor the Rule of Necessitation.\footnote{That is, the system $\JLLuk$ is a nonnormal modal logic by design.} These lacks are by design, but they are not \textit{ad hoc} since they obey the goal of accounting for epistemic processes of non-idealized agents with bounded cognitive resources, who do not posses knowledge of all the tautologies of propositional calculus. Most of the time, if these agents want to know a tautology, they have to deduce it, and even if it's a tautology, it would not count as knowledge until it has been properly proven or justified. The system $\JLLuk$ captures these facts. We characterize the system $\JLLuk$ in the following section.

\subsection{System $\JLLuk$} 
The majority of justification logic systems take the system $\mathsf{J_0}$ as the point of departure (which is axiomatically defined in \cite{ArtemovFitting2019}), but since $\mathsf{J_0}$ includes all the tautologies of the language of classical propositional logic among its axioms, we are going to deviate from this point of departure because the goal of our system is to derive tautologies and theorems, showing their epistemic contribution. However, like most justification logic systems, we shall include \textit{Application} and \textit{Sum} axioms, plus \textit{Modus Ponens} as an inference rule. Since a justification logic system equipped with nothing more than \textit{Application}, \textit{Sum} and \textit{Modus Ponens} would be even weaker than $\mathsf{J_0}$, which cannot prove that any formula has a justification,\cite{artemov2005introducing} we will extend it with schemes for \L ukasiewicz's axioms accompanied by their respective justifications ---hence the name ``System $\JLLuk$''. We have chosen \L ukasiewicz's axioms because this set of axiom schemes has enough power to derive the set of all tautologies of propositional logic.\footnote{The choice of \L ukasiewicz's axioms may affect the length of the derivations, but since there is no polynomially bounded propositional proof system, we hypothesize that our system is a good as any other one in terms of the length of derivations.} The standard definitions of a justification logic system (justification terms, \textit{Application} and \textit{Sum} axioms, constant specification, etc.) are borrowed from \cite{justification:logic,artemov2005introducing,ArtemovFitting2019}.     

\begin{defn}[Justification terms] The set $T$ of justification terms is built from countably many constants  in the set $C$ of constant justification terms, and countably many variables in the set $V$ of variable justification terms, given by
\[ 
t ::= c \mid x \mid [t \cdot t] \mid [t + t]
\]
where $c \in C$ and $x \in V$.\footnote{Following \cite{ArtemovFitting2019}, we will use brackets ``[]'' for terms, and parenthesis ``()'' for formulas.}  
\end{defn}

\begin{defn}Justification formulas] Justification formulas $\phi, \psi$ in the set of $F$ of justification formulas are built according to the following grammar
\[ 
\phi, \psi ::= p \mid \neg \phi \mid \phi \lor \psi \mid  \phi \land \psi \mid \phi \to \psi \mid t:\phi
\]
where $p$ is in the set $P$ of countably many atomic propositions, connectives are defined in the standard way, and $t$ is in the set $T$. The justification formula $t:\phi$ is read as `the terms $t$ is a justification for $\phi$.' 

\end{defn}

\begin{defn}[System $\JLLuk$] Let $\phi, \psi$ be formulas  in the set of $F$ of justification formulas, and $t, s$ terms in the set $T$ of justification terms. The system is defined by the \textit{Application} and \textit{Sum} axioms, schemes for \L ukasiewicz's axioms with a justification constant for each instance of these axioms, \textit{Modus Ponens} (MP) as an inference rule, plus one additional rule to infer instances of justified applications, which we call \textit{Application introduction}:\footnote{\textit{Application introduction} is an inferential rule we derived from the \textit{Application} axiom  justification logic; this new rule is equivalent to call the axiom and then apply  \textit{Modus Ponens}.}
\begin{align*}
    & \textit{Application}  && s:(\phi \to \psi) \to (t:\phi \to [s\cdot t]:\psi)\\
    & \textit{Sum}  && s:\psi \to [s + t]:\psi, t:\psi \to [s + t]:\psi\\
    &\text{\L ukasiewicz's \ Axiom 1} &&c_1: \phi \to (\psi \to \phi)\\
    &\text{\L ukasiewicz's \ Axiom 2} && c_2: \phi \to (\psi \to \chi) \to ((\phi \to \psi) \to (\phi \to \chi))\\
    &\text{\L ukasiewicz's \ Axiom 3} && c_3: (\neg \phi \to \neg \psi) \to (\psi \to \phi)\\
    & \textit{MP}\ && \phi \to \psi,\ \phi  \vdash \psi\\
    & \textit{Application introduction}\ &&  s:(\phi \to \psi),\  t:\phi,\ \vdash  [s\cdot t]:\psi\\
\end{align*}

\end{defn}

\L ukasiewicz's axiom schemes can be instantiated as many times as needed. Note that for each instance of these axiom schemes, there is a constant specification $c$ that justifies it. Regarding the condition on the set of constant specifications, we are going to work with a schematic set, that is, if $A$ and $B$ are both instances of the same axiom scheme, $c:\ A \in CS$ if and only if $c:\ B \in CS$, for every constant symbol $c$. Hence, the system $\JLLuk$ will have a constant specification schematic numerable set of formulas $CS$ to  assign a unique constant for every axiom. The working of the system $\JLLuk$ is specified as follows:

\begin{defn}[Consequence] Let $\JLLuk$ be a justification logic system with $CS$ as its constant specification, and an arbitrary set of formulas $\Gamma$. Then $\Gamma  \vdash_{\JLLuk} \phi$ if there is a finite sequence of formulas ending with $\phi$, in which each formula is either an axiom of $\JLLuk$, a member of $CS$ or $\Gamma$, or follows from earlier formulas by applying the inferential rules allowed by the system. The system benefits from the deduction theorem, that is, for any $\Gamma \subseteq T$ and $\phi, \psi \in T$, it is true that $\Gamma, \phi \vdash_{\JLLuk} \psi$ if and only if  $\Gamma \vdash_{\JLLuk} \phi \to \psi$. 

\end{defn}

The system $\JLLuk$ is the result of enriching a fragment of $\mathsf{J_0}$ ---that which does not include all tautologies of propositional logic among its axioms--- with \L ukasiewicz's axioms, and since these axioms are valid, $\JLLuk$ is sound and complete with respect to \L ukasiewicz's  proof system,\cite{Lukasiewicz,Lukasiewicz:selected} which in turn is sound and complete for the language of propositional logic.\footnote{It is well known that everything that can be proven using \L ukasiewicz's proof system within the language of propositional logic is valid in this language, and everything that is valid in this language can be proven by means of \L ukasiewicz's proof system.} It will be easy to see that the system $\JLLuk$ preserves the soundness and completeness from \L ukasiewicz's  proof system: if something is provable by means of instances of the \textit{Application} or \textit{Sum} axioms and justified instances of  \L ukasiewicz's axioms, then it is valid in $\JLLuk$ and if something is valid in $\JLLuk$, then it is provable by means of instances of the \textit{Application} or \textit{Sum} axioms and justified instances of  \L ukasiewicz's axioms (see theorems \ref{theo:correcteness}, \ref{theo:correcteness:extended}).

Hence, for every theorem that is a justified formula in $\JLLuk$, the formula without the justification is a theorem in \L ukasiewicz's  proof system, that is, if $t: \alpha$ is a theorem in $\JLLuk$, then $\alpha$ is a theorem in \L ukasiewicz's  proof system. This claim is illustrated through the two following derivations in which $q \to r$ is deduced from the premises $p, q \to (p\to r)$; the first of this derivations is in $\vdash_{\Lslash}$ and the second one in $\vdash_{\JLLuk}$:

\begin{align*}
  & \textbf{Derivation in $\vdash_{\Lslash}$}\\ 
  & 1.\ p   & & \text{premise}\\
  & 2.\ q \to (p \to r)   & & \text{premise}\\
  & 3.\ p \to (q \to p)   & & \text{\L ukasiewicz's \ Axiom 1}\\
  & 4.\ (q \to p) & & \text{MP 1, 3}\\
  & 5.\ (q \to (p \to r)) \to ((q \to p) \to (q \to r)) & & \text{\L ukasiewicz's \ Axiom 2}\\
  & 6.\ (q \to p) \to (q \to r)  & & \text{MP 2, 5}\\
  & 7.\ q \to r  & & \text{MP 4, 6}
\end{align*}

\begin{align*}
  & \textbf{Derivation in $\vdash_{\JLLuk}$}\\
  & 1.\ x:p   & & \text{premise}\\
  & 2.\ y:q \to (p \to r)   & & \text{premise}\\
  & 3.\ c_1:p \to (q \to p)   & & \text{$\JLLuk$ Axiom 1}\\
  & 4.\ [c_1 \cdot x]:(q \to p) & & \text{Application introduction 3,1}\\
  & 5.\ c_2:(q \to (p \to r)) \to\\
  &\qquad\quad ((q \to p) \to (q \to r)) & & \text{$\JLLuk$ Axiom 2}\\
  & 6.\ [c_2 \cdot y]:(q \to p) \to (q \to r)  & & \text{Application introduction 5,2}\\
  & 7.\ [c_2 \cdot y]\cdot[c_1 \cdot x]:q \to r  & & \text{Application introduction 4,6}
\end{align*}

The deduction in $\vdash_{\Lslash}$ is straightforward: $q \to r$ is deduced from $p$ and $q \to (p \to r)$ using instances of \L ukasiewicz's axiom schemes and  \textit{Modus Ponens}. The deduction in $\vdash_{\JLLuk}$ illustrates how to deduce $q \to r$ explicitly showing the justification of accepting $q \to r$. In this proof, all the tools available in the set $T$ of justification terms from definition 3.1 above are used: variables are used to justify the given premises, while constant specifications are used to justify instances of \L ukasiewicz's axioms. Then, to prove that $ q \to r$ is justified in $\JLLuk$, we rely on the \textit{Application} axiom of justification logic and the \textit{Application introduction} rule to derive a justification for the formula $ q \to r$; in this way we derived $[c_2 \cdot y]\cdot[c_1 \cdot s]:q \to r$. From this, we can generalize that a formula $\phi$ that can be deduced in $\vdash_{\Lslash}$ if and only if there is a justification term $t$ such that $\vdash_{\JLLuk} t: \phi$; this result will be proved in theorem 3.1 bellow.  

As the examples above show, it is easy to transform a  $\vdash_{\Lslash}$ proof into a  $\vdash_{\JLLuk}$ proof. The following theorem generalizes this idea.  

\begin{theorem}\label{theo:correcteness}
$\vdash_{\Lslash} \beta$ if and only if there exists a justification term $t$ such that
$\vdash_{\JLLuk} t:\beta.$
\end{theorem}

\begin{proof}
From $\vdash_{\Lslash } \beta$ we can prove $\vdash_{\JLLuk} t:\beta$ by induction on proofs in $\vdash_{\Lslash}$:
\begin{itemize}
    \item[(i)] if $\beta$ is an instance of an axiom, by definition we have a specification constant $c$ such that $c:\beta$;
    \item[(ii)] if $\beta$ was derived by MP from $\alpha \to \beta$ and $\alpha$, then by inductive hypothesis there are terms $s$ and $u$ such that $s: \alpha \to \beta$ and $u:\alpha$. Then by the Application axiom and two applications of MP we have $[s\cdot u] : \beta$.
\end{itemize}

On the other hand, if $\vdash_{\JLLuk} t:\beta$, we can prove $\vdash_{\Lslash} \beta$ by structural induction on $t$, as there are three options:

\begin{itemize}
    \item[(a)] $t = c$ in which case $\beta$ is an (instance of an) axiom and clearly $\vdash_{\Lslash} \beta$;
    \item[(b)] $t = [s \cdot u]$ and so previously we had $\vdash_{\JLLuk}s: \alpha \to \beta$ and $\vdash_{\JLLuk}u: \alpha$ for some $\alpha$ and, by induction hypothesis, $\vdash_{\Lslash} \alpha \to \beta$ and $\vdash_{\Lslash} \alpha$. Then, by \textit{Modus Ponens} we have $\vdash_{\Lslash} \beta$;
    \item[(c)] $t = [s+u]$ in which case we previously had $s : \beta$ or $u: \beta$, by induction hypothesis then $\vdash_{\Lslash} \beta$ in either case;
    \item[(d)] there is really no case $t = x$, as variable terms do not appear in proofs without premises.

\end{itemize}
\end{proof}

Observe how the ``if'' direction of this theorem is analogous to the \emph{internalization property} of traditional justification logics \cite[20]{ArtemovFitting2019}, while the ``only if'' direction is analogous to \emph{inversed internalization} \cite[108]{goetschi}.

\begin{theorem}\label{theo:correcteness:extended}
$$\alpha_1, \dots, \alpha_n \vdash_{\Lslash} \beta$$
if and only if there exist a justification term $t$ such that
$$x_{i_1}: \alpha_1, \dots, x_{i_n}: \alpha_n  \vdash_{\JLLuk} t:\beta.$$
\end{theorem}

\begin{proof}
It is very similar to the previous case, but because we are adding premises that may not be (instances of) axioms, we introduce variable justification terms. In one direction we need to introduce a (iii) case where $\beta$ is a premise (and therefore, its justification term would be a variable). In the opposite direction, we have to consider a (d) case, because now it is possible that $t = x$, in which case $\beta$ is a premise and therefore it is trivial that $\beta \vdash_{\Lslash} \beta$.
\end{proof}

Again, this theorem is analogous to the Lifting Lemma \cite[22]{ArtemovFitting2019} and inversed internalization \cite[112]{goetschi}.

The two proofs in $\vdash_{\Lslash}$ and $\vdash_{\JLLuk}$ above exemplify how a formula that was proven in $\vdash_{\Lslash}$ can be proven in $\vdash_{\JLLuk}$ as well, but with the advantage that the proof in $\vdash_{\JLLuk}$ makes explicit the reasons that justify the formula. In this way, the system $\JLLuk$ exhibits the epistemic contribution of the deduction or derivation of a formula.

\section{Conclusions and future work}

In this paper we have developed an alternative account of the epistemic contribution of deduction by exhibiting the process of deducing a theorem of the language of propositional logic, which in turn, allows the subject to know the deduced theorem. To this end, we have addressed the following challenges:

To show that the process of deduction bears an epistemic contribution, in Section 3 we developed a justification logic system extended with schemes for \L ukasiewicz's axioms: the system $\JLLuk$. The system $\JLLuk$ allowed us to explicitly express within the language the reasons or justification for accepting a theorem, which in turn reflects the process that leads to knowledge acquisition through deduction. 
    
In this way, we have provided an alternative account for the epistemic contribution of deduction that sheds new light upon classical problems within the philosophy of logic such as Hintikka's scandal of deduction and Cohen and Nagel's paradox of inference. Moreover, our proposal not only provides theoretical tools to exhibit and measure the epistemic contribution of deduction, but it also improves our understanding of deduction as a process of providing reasons to justify the logical consequences of our current knowledge. Therefore, the insights gained from this work expand our understanding of knowledge by deduction and its potential as a source of new knowledge.

Future work should focus on the complexity of deciding logical entailment and its implications for the epistemic contribution of deduction. It is known that, in terms of computational complexity, deciding whether $\varGamma \vdash \alpha$ is not a trivial task since it is at least as hard as deciding whether a formula of propositional logic is a tautology, which is a decision problem known as {\sc taut} and belongs to the co-{\sc np} complexity class. Hence, there is no efficient algorithm that can always tell whether a formula $\alpha$ is implied by a $\varGamma$. Furthermore, even if it is the case that $\varGamma \vdash \alpha$, we may not be able to verify it in polynomial time. This shows another way in which getting to know the logical consequences of known premises is complex enough. In addition, further research should also focus on the complexity of finding a justification term $t$ such that $\varGamma \vdash_{\JLLuk} t:\alpha$. Moreover, finding the shortest term $t$ such that $\varGamma \vdash_{\JLLuk}t:\alpha$ is even harder, as at the very least it implies solving {\sc taut}, but very likely more as right now it is known to be an {\sc np}-hard problem \cite{minimal-proof-length,cook1971complexity}.

Being limited to knowledge by deduction, this work has not covered any other form of inference or argument, though the main results may be generalizable. Notwithstanding this limitation, this work expands the repertoire of epistemic formal tools by presenting a novel system $\JLLuk$ to account for the epistemic contribution of deduction. To the knowledge of the authors, formal tools from justification logic have never been used to provide an account of the epistemic contribution of deduction. In this way, the present work provides a useful framework to develop an unprecedented and robust understanding of deduction and its epistemic contribution.  



\bibliographystyle{apalike}
\bibliography{references.bib}

\begin{thebibliography}{}

\bibitem[Abramsky, 2008]{ABRAMSKY2008}
Abramsky, S. (2008).
\newblock Information, processes and games.
\newblock In Adriaans, P. and {van Benthem}, J., editors, {\em Philosophy of
  Information}, Handbook of the Philosophy of Science, pages 483--549.
  North-Holland, Amsterdam.

\bibitem[Alekhnovich et~al., 2001]{minimal-proof-length}
Alekhnovich, M., Buss, S., Moran, S., and Pitassi, T. (2001).
\newblock Minimum propositional proof length is np-hard to linearly
  approximate.
\newblock {\em The Journal of Symbolic Logic}, 66(1):171--191.

\bibitem[Arora and Barak, 2009]{arora_barak_2009}
Arora, S. and Barak, B. (2009).
\newblock {\em Computational Complexity: A Modern Approach}.
\newblock Cambridge University Press.

\bibitem[Artemov, 1994]{artemov1994}
Artemov, S. (1994).
\newblock Logic of proofs.
\newblock {\em Annals of Pure and Applied Logic}, 67(1-3):29--59.

\bibitem[Artemov, 1995]{artemov1995}
Artemov, S. (1995).
\newblock Operational modal logic.
\newblock Technical Report~29, Mathematical Sciences Institute, Cornell
  University.

\bibitem[Artemov, 2008]{justification:logic}
Artemov, S. (2008).
\newblock The logic of justification.
\newblock {\em The Review of Symbolic Logic}, page 477–513.

\bibitem[Artemov and Fitting, 2019]{ArtemovFitting2019}
Artemov, S. and Fitting, M. (2019).
\newblock {\em Justification Logic: Reasoning with Reasons}, volume 216.
\newblock Cambridge University Press.

\bibitem[Artemov and Kuznets, 2009]{ArtemovKuznets2009}
Artemov, S. and Kuznets, R. (2009).
\newblock Logical omniscience as a computational complexity problem.
\newblock In {\em Proceedings of the 12th Conference on Theoretical Aspects of
  Rationality and Knowledge}, pages 14--23. TARK.

\bibitem[Artemov and Nogina, 2005]{artemov2005introducing}
Artemov, S. and Nogina, E. (2005).
\newblock Introducing justification into epistemic logic.
\newblock {\em Journal of Logic and Computation}, 15(6):1059--1073.

\bibitem[Carnap and Bar-Hillel, 1952]{Carnap1952CSI}
Carnap, R. and Bar-Hillel, Y. (1952).
\newblock An outline of a theory of semantic information.
\newblock Technical Report 247, Research Laboratory of Electronics,
  Massachusetts Institute of Technology.

\bibitem[Cohen and Nagel, 1934]{Cohen&Nagel}
Cohen, M. and Nagel, E. (1934).
\newblock {\em An introduction to logic and scientific method}.
\newblock Routledge and Kegan Paul.

\bibitem[Cook, 1971]{cook1971complexity}
Cook, S.~A. (1971).
\newblock The complexity of theorem-proving procedures.
\newblock In {\em Proceedings of the third annual ACM symposium on Theory of
  computing}, pages 151--158.

\bibitem[Dummett, 1978]{Dummett1978}
Dummett, M. A.~E. (1978).
\newblock {\em Truth and Other Enigmas}.
\newblock Harvard University Press.

\bibitem[Du{\v{z}}{í}, 2010]{Duvzi2010paradox}
Du{\v{z}}{í}, M. (2010).
\newblock The paradox of inference and the non-triviality of analytic
  information.
\newblock {\em Journal of Philosophical Logic}, 39(5):473--510.

\bibitem[D’Agostino and Floridi, 2009]{floridi:dagostino:scandal}
D’Agostino, M. and Floridi, L. (2009).
\newblock The enduring scandal of deduction: Is propositional logic really
  uninformative?
\newblock {\em Synthese}, 167:271--315.

\bibitem[Goetschi, 2012]{goetschi}
Goetschi, R. (2012).
\newblock {\em On the Realization and Classification of Justification Logics}.
\newblock PhD thesis, University of Bern.

\bibitem[Hintikka, 1970]{hintikka1970}
Hintikka, J. (1970).
\newblock Information, deduction, and the a priori.
\newblock {\em Nous}, pages 135--152.

\bibitem[\L{}ukasiewicz, 1964]{Lukasiewicz}
\L{}ukasiewicz, J. (1964).
\newblock {\em Elements of Mathematical Logic}.
\newblock New York: Macmillan.

\bibitem[{\L ukasiewicz}, 1970]{Lukasiewicz:selected}
{\L ukasiewicz}, J. (1970).
\newblock {\em Selected works}.
\newblock North Holland Publishing Company.

\bibitem[Mill, 1843]{Mill1843}
Mill, J.~S. (1843).
\newblock {\em Mill, A System of Logic, Ratiocinative and Inductive, Being a
  Connected View of the Principles of Evidence and the Methods of Scientific
  Investigation}.
\newblock Longmans, Green, and Co.

\bibitem[Parikh, 2008]{parikh2008}
Parikh, R. (2008).
\newblock Sentences, belief and logical omniscience, or what does deduction
  tell us?
\newblock {\em The Review of Symbolic Logic}, 1(4):459--476.

\bibitem[Primiero, 2007]{primiero2007information}
Primiero, G. (2007).
\newblock {\em Information and knowledge: a constructive type-theoretical
  approach}, volume~10.
\newblock Springer Science \& Business Media.

\bibitem[Rumfitt, 2008]{Rumfitt2008}
Rumfitt, I. (2008).
\newblock Knowledge by deduction.
\newblock {\em Grazer Philosophische Studien}, 77(1).

\bibitem[Sequoiah-Grayson, 2008]{sequoiah2008}
Sequoiah-Grayson, S. (2008).
\newblock The scandal of deduction.
\newblock {\em Journal of Philosophical Logic}, 37(1):67--94.

\bibitem[van Benthem and Martinez, 2008]{BenthenMartinez}
van Benthem, J. and Martinez, M. (2008).
\newblock The stories of logic and information.
\newblock In Adriaans, P. and {van Benthem}, J., editors, {\em Philosophy of
  Information}, Handbook of the Philosophy of Science, pages 217--280.
  North-Holland, Amsterdam.

\end{thebibliography}

\end{document}